\documentclass[11pt,a4paper,leqno]{article}
\usepackage{amsfonts,amssymb,mathrsfs,amsmath}

\usepackage{amscd}

\vspace{0.5cm}

 4

\font\erm=cmr8

\usepackage{graphicx}

\author{Andrzej Krzysztof Kwa\'sniewski}
\title{A note on Ward-Horadam \textit{H(x)}-binomials' recurrences}

\newtheorem{defn}{Definition}

\newtheorem{theoremn}{Theorem}

\newtheorem{quot}{Quotation}

\newcommand{\fnomial}[2]{ {{#1} \choose {#2}}_{\!\!F} }
\newcommand{\fnomialF}[3]{ {{#1} \choose {#2}}_{\!\!#3} }

\begin{document}

\begin{center}
\noindent {\sffamily A note on Ward-Horadam \textit{H(x)}-binomials' recurrences}  \\ 

\vspace{0.5cm} 
\sffamily Andrzej Krzysztof Kwa\'sniewski

\vspace{0.5cm}
{\erm
\sffamily Member of the Institute of Combinatorics and its Applications, 
Winnipeg, Manitoba, Canada \\
	PL-15-674 Bia\l ystok, Konwaliowa 11/11, Poland\\
	e-mail: kwandr@gmail.com\\}
\end{center}

\vspace{0,5cm}

\noindent \textbf{Summary} 

\noindent \sffamily As a matter of continuation of  \cite[2010]{akk2010} - we deliver here  \textit{H(x)}-binomials' recurrence formula appointed by Ward-Horadam $ H(x)= \left\langle H_n(x)\right\rangle _{n\geq 0}$ functions' sequence which comprises   in $H \equiv H(x=1)$ number sequences case the  \textit{V}-binomials' recurrence formula determined by the primordial  Lucas sequence of the second kind   $V = \left\langle V_n\right\rangle_{n\geq 0}$ as well  as its well elaborated companion fundamental Lucas sequence of the first kind $ U = \left\langle U_n\right\rangle_{n\geq 0}$ which gives rise in its turn to the  \textit{U}-binomials' recurrence as in \cite[1878]{EdL} , \cite[1949]{JM},  \cite[1964]{TF},  \cite[1969]{Gould},  \cite[1989]{G-V}  or in  \cite[1989]{K-W}  etc. 

\vspace{0,3cm} 

\noindent For the sake of combinatorial interpretations and in number theory  $H(x=1)$, $H_n(x=1) \equiv H_n$ is  usually considered to be natural or integer numbers valued sequence. Number sequences $H =H(x=1)= \left\langle H_n\right\rangle _{n\geq 0}$ were recently called by several authors: Horadam sequences. 

\vspace{0,2cm}

\noindent The list of references is mostly indicatory (see references therein) and is far from being complete.

\vspace{0,3cm}
\noindent AMS Classification Numbers: 05A10 , 05A30.

\vspace{0,2cm}

\vspace{0.5cm}


\section{General Introduction}

\noindent \textbf{1.1.  $p,q$ people are Lucas' followers people.} The are many authors who use in their investigation the fundamental Lucas sequence $ U \equiv \left\langle n_{p,q}\right\rangle_{n\geq 0}$ - frequently  with different notations -  where $n_{p,q} = \sum_{j=0}^{n-1}{p^{n-j-1}q^j}=U_n$; see Definition 1  and then definitions that follow it.  In regard to this a brief intimation   is on the way.

\vspace{0,1cm}

\noindent   To our knowledge it was  Fran\c{c}ois  \'Edouard Anatole Lucas in \cite[1878]{EdL} who was the first who had  \textit{not only} defined \textit{fibonomial} coefficients as stated in \cite[1989]{K-W}   by Donald  Ervin  Knuth and Herbert Saul Wilf but who was the first who had defined  \textit{$U_n \equiv n_{p,q}$}-binomial coefficients ${n \choose k}_U \equiv {n \choose k}_{p,q}$  and had derived a recurrence for them: see page 27, formula (58) \cite[1878]{EdL}. Then - referring to Lucas -  the investigation relative to  divisibility properties of relevant number Lucas sequences $D$, $S$ as well as numbers' $D$ - binomials and  numbers' $D$ - multinomials was continued in \cite[1913]{Carmichel}  by  Robert Daniel Carmichel;  see pp. 30, 35 and 40 in \cite[1913]{Carmichel}  for  $ U \equiv D = \left\langle D_n\right\rangle_{n\geq 0}$ and ${n \choose {k_1,k_2,...,k_s}}_D$  - respectively. Note there also formulas (10) , (11) and (13)  which might perhaps  serve to derive explicit untangled form of recurrence for the  $ \textit{V}$- \textit{binomial coefficients} ${n \choose k}_V \equiv {n \choose k}_S$ denoted by  primordial Lucas sequence $\left\langle S_n\right\rangle_{n\geq 0}= S \equiv V$. Number sequence $F(x=1 = A)$  $A$ - multinomial coefficients' \textit{recurrences} are not present in that early works and up  to our knowledge  
 a special case of such appeared at first in \cite[1979]{Shannon 1979 multi} by  Anthony G. Shannon. More on that - in what follows after Definition 3.

\vspace{0,1cm}

\noindent Significant peculiarity of Lucas originated sequences  includes  their importance for number theory  (see  middle-century paper \cite{Halton 1996} by  John H. Halton and recent, this century papers \cite[2010]{Smyth1} by Chris Smith  and \cite[2010]{Smyth2} by K\'alm\'an  Gy{\"{o}}ry with Chris Smith  and the reader may enjoy also the PhD Thesis  \cite[1999]{Anne-Marie 1999}  by  Anne-Marie Decaillot-Laulagnet). This Lucas originated investigation amalgamates diverse areas of mathematics due to hyberbolic - trigonometric character of these Fonctions Num\'eriques Simplement Périodiques i.e. fundamental and primordial Lucas sequences - as beheld  in \cite[1878]{EdL}. One may track then  a piece of further constructions  for example in \cite[1999]{akk1999}).

\noindent There in \cite[1999]{akk1999} tail formulas (3.12) and (3.14)  illustrating proved and exploited by \'Eduard Lucas  the complete analogy of the $V_n$  and $U_n$  symmetric functions of roots with the circular and hyperbolic functions of order $ 2$ due to Lucas formulas (5) in \cite{EdL} rewritten in terms of $\cosh$  and $\sinh$ functions as formulas (3.13) and (3.14) in \cite{akk1999} as resulting from   de Moivre one parameter group  introduced in \cite{akk1999} via (1.4) in order to pack compactly  addition formulas (1.6), (1.7) in \cite{akk1999}  equivalent to  (49)  and corresponding  recurrence relations in \cite{EdL} into abelian group  "parcel"   encompassing  Tchebycheff polynomials of both kinds. 

\vspace{0,1cm}

\noindent In this connection see the Section 2 in the  recent Ward-Horadam people paper  \cite [2009]{He-Siue}  by  Tian-Xiao He,  Peter Jau-Shyong Shiue. There in Proposition 2.7.  illustrative   Example 2.8. with Tchebycheff polynomials of the first kind the well known recurrence formula (2.28)  is equivalent to  abelian one-parameter de Moivre matrix group  multiplication rule from which 
the corresponding recurrence  (1.7)   in  \cite[1999]{akk1999}  follows.

\vspace{0,3cm}


\noindent \textbf{2.1.  As has been foreshadowed} in \cite[2010]{akk2010} we deliver here - continuing  the note \cite{akk2010} - the \textit{H(x)}-binomials' recurrence formula appointed by Ward-Horadam $ H(x)= \left\langle H_n(x)\right\rangle _{n\geq 0}$ field of zero characteristic nonzero valued functions' sequence which comprises  for $H \equiv H(x=1)$ number sequences case - the  \textit{V}-binomials' recurrence formula determined by the primordial Lucas sequence of the second kind  $V = \left\langle V_n\right\rangle_{n\geq 0}$   \cite[2010]{akk2010} as well  as its well elaborated companion fundamental Lucas sequence of the first kind $ U = \left\langle U_n\right\rangle_{n\geq 0}$ which gives rise in its turn to the  \textit{U}-binomials' recurrence as in \cite[1878]{EdL} , \cite[1949]{JM},  \cite[1964]{TF},  \cite[1969]{Gould},  \cite[1989]{G-V}  or in  \cite[1989]{K-W}  and so on. 

\vspace{0,3cm}

\noindent We do it by following recent applicable work \cite[2009]{Savage}  by Nicolas A. Loehr and Carla  D.  Savage   thought one may  - for that purpose - envisage now easy extensions of particular $p,q$ - cases considered earlier - as  \textit{for example} the following: the relevant recursions in \cite[1949]{JM}, in  \cite[1989]{G-V},  in \cite[1992]{KaKi} - ( recursions  (40) and (51)) , or \cite[2000]{Holte 2000} by  John  M.  Holte (Lemmas 1,2 dealing with  $U$-binomials provide a  motivated example for observation  Theorem 17 in \cite{Savage} ) One is invited also to track Lemma 1 in  \cite[2001]{Hu-Sun} by  Hong Hu  and   Zhi-Wei Sun ; see also corresponding recurrences for $p,q$-binomials $\equiv$ $U$-binomials in \cite[1878]{EdL} or in \cite[2008]{MD2} \textbf{v[1]} by Maciej Dziemia\'nczuk (compare there (1) and (2) formulas), or see Theorem 1  in \cite[2008]{Corsino}  by Roberto Bagsarsa Corcino as well as  track the proof of the Corollary 3. in \cite[2009]{MD2} \textbf{v[2]}  by Maciej Dziemia\'nczuk.

\vspace{0,3cm}



\noindent \textbf{This looked for }here \textit{H(x)}-binomials' recurrence formula (recall: encompassing \textit{V}-binomials for primordial Lucas sequence $V$) is not present neither in \cite{EdL} nor in \cite{Savage}, nor in \cite[1915]{Fon} , nor in \cite[1936]{Ward}, nor in \cite[1949]{JM}. Neither we find it in  - quoted here contractually by a nickname as  "Lucas $(p,q)$-people" - references [1-44]. Neither it is present in  all other  - quoted here contractually by a nickname as  "Ward-Horadam -people" - references [49-79]. Ad "Lucas $(p,q)$-people" and "Ward-Horadam -people" references -
(including these [n]  with  $n>73$- the distinction  which are which is  quite contractual. The nicknames are nevertheless indicatively helpful.  We shall be more precise soon - right with definitions are being started.



\vspace{0,3cm}

\noindent \textbf{Meanwhile}  \textit{H(x)}-binomials' recurrence formula for the  Ward-Horadam sequence  $H(x) = \left\langle H_n(x)\right\rangle _{n\geq 0}$  follows straightforwardly from the easily proved important observation - the Theorem 17 in \cite[2009]{Savage} as already had it been remarked in \cite[2010]{akk2010} for the $H \equiv H(x=1)$ case. 

\vspace{0,3cm}

\noindent \textbf{This paper formula} may and should be confronted with  Fonten\'e obvious recurrence for complex valued $A$-binomials ${n \choose k}_A$,  $A \equiv A(x=1)$  in \cite[1915]{Fon} i.e. with (6) or (7) identities in \cite[1969]{Gould} by Henri W. Gould or with recurrence in \cite[1999]{ALs} by Alexandru Ioan  Lupas ,  which particularly also stem easily just  from the definition of \textbf{any} $F(x)$-binomial coefficients arrays with $F(x) = \left\langle F_n(x)\right\rangle_{n\geq 0}$ staying for any field of characteristic zero  nonzero valued functions'sequence ; $F_n(x) \neq 0,\:n\geq 0$.  For $F = F(x=1)$-multinomial coefficients automatic definition see \cite[1913]{Carmichel}  by    Robert Daniel Carmichel or then  \cite[1969]{Gould} by Henri W. Gold and finally  see \cite[1979]{Shannon 1979 multi} by Anthony G. Shannon, where recurrence is proved for ${n \choose {k_1,k_2,...,k_s}}_U$    with $U$-Lucas fundamental being here complex valued number sequence. For $F(x)$ - multinomial coefficients see  \cite[2004]{akk3 2004}  and compare with  $F(x)$-binomials from \cite[1999]{ALs} or those from \cite[2001]{Rich 2001}.

\vspace{0,2cm}

\noindent To this end we supply now two informations pertinent  ad references and ad nomenclature.

\vspace{0,3cm}

\noindent \textbf{3.1. Ad the number theory and divisibility properties references.}  
\noindent For the sake of combinatorial interpretations of $F$- number sequences as well as their correspondent $F$-multinomial coefficients  and also for the sake of the number theoretic studies of Charles Hermite  \cite{Herm 1889-1} and with Thomas Jan Stieltjes in   \cite{Herm 1889-2} or by  Robert Daniel Carmichel  \cite[1913]{Carmichel} or  \cite[1919]{Dick} or that of Ward  \cite[1936]{Ward1}, \cite[1939]{Ward2}, \cite[1937]{Ward3}, \cite[1937]{Ward4}, \cite[1954]{Ward0}, \cite[1955]{Ward5}, \cite[1959]{Ward6} and that of  Lehmer  \cite[1930]{Leh}, \cite[1933]{Leh1},  \cite[1935]{Leh2} or this of  Andrzej Bobola Maria Schinzel   \cite[1974]{Schinzel} and Others' studies  \textit{on  divisibility properties} -  these are the sub-cases  $F_n \in\mathbb{N}$  or $F_n \in\mathbb{Z}$ which are being regularly considered at the purpose.  

\vspace{0,1cm}

\noindent   As for the "Others" - see \textit{for example}:  \cite[1959]{Durst}, \cite[1973]{RAKKK},  \cite[1974]{SingMa1},  \cite[1974]{SingMa2},  \cite[1974]{SingMa3},  \cite[1973]{SingMa4}, \cite[1977]{KKK},  \cite[1977]{Stewart}, \cite[1979]{ShaHor},  \cite[1979]{HorLohSha 1979},  \cite[1980]{Somer},  \cite[1980]{Beukers},  \cite[1989]{K-W}, \cite[1991]{SatoShi}, \cite[1992]{Zhi-Wei Sun}, \cite[1995]{GouldPaula},  \cite[1999]{adGould},  \cite{GouldList},   \cite[1995]{Hil-Vranc}, \cite[1995]{KimWeb}, \cite[1995]{Voutier},  \cite[1998]{Wilson},    \cite[2001]{Hu-Sun}, \cite[2006]{MoAbo}, \cite[2009]{MDWB}. 

\vspace{0.3cm}

\noindent \textbf{3.4. Ad Ward-Horadam naming.} 
\noindent According to the authors of  \cite[2009]{He-Siue}  it was  Mansour   \cite{Mansour}  who called  the  sequence $H = \left\langle a_n \right\rangle n\geq 0$  defined by (1)  a \textbf{Horadam's sequence}, as - accordingly to the author of \cite{Mansour} - the number sequence $H$ was introduced in 1965 by Horadam  \cite{Horadam1} (for special case of Ward-Horadam number sequences see Section 2 in \cite[1974]{WaltHorad} and see also  \cite[2009]{HorzumKocer}), this however  notwithstanding  the ingress of complex numbers valued $F$-\textit{binomials}  and $F$-\textit{multinomials} into Morgan Ward's systematic \textit{Calculus of sequences} in  \cite[1936]{Ward} and then in 1954 Ward's introduction   of "`\textit{nomen omen}"'  $W\equiv H$ in \cite[1954]{Ward0} integer valued sequences. 

\vspace{0,1cm}

\noindent Perceive then the appraisal of adequate Morgan Wards' work in the domain by  Henri W.Gould \cite[1959]{Gould} and by Alwyn F. Horadam   and Anthony G. Shannon  in  \cite[1976]{HorSha}  or Derrick  Henry Lehmer in \cite[1993]{Leh3}. On this occasion note also the \textbf{Ward-Horadam} number sequences \textbf{in} \cite[1965]{Zeitlin0}  and  \cite[1965]{Zeitlin}.

\section{Preliminaries}

\vspace{0.1cm}

\textit{Names}:  The Lucas sequence   $V = \left\langle V_n\right\rangle_{n\geq 0}$  is called the Lucas sequence of the second kind - see: \cite[1977, Part I]{KKK}, or \textbf{primordial} - see
\cite[1979]{ShaHor}.     

\noindent The Lucas sequence  $ U = \left\langle U_n\right\rangle_{n\geq 0}$   is called the Lucas sequence of the first kind - see: \cite[1977, Part I]{KKK}, or \textbf{fundamental}  - see p. 38 in \cite[1949]{JM} or see \cite[1979]{Shannon 1979 multi}  and  \cite[1979]{ShaHor}.     

\vspace{0.2cm} 

\noindent In the sequel we shall deliver the looked for recurrence for $H$-binomial coefficients ${n \choose k}_H$  determined by the Ward-Horadam sequence $H$ - defined below.

\vspace{0.2cm}

\noindent In compliance with Edouard Lucas' \cite[1878]{EdL} and twenty, twenty first century $p,q$-people's notation we shall at first review here in brief the general second order recurrence;
(compare this review  with the  recent "Ward-Horadam" peoples' paper  \cite [2009]{He-Siue} by Tian-Xiao He and Peter Jau-Shyong Shiue or earlier  $p,q$-papers \cite[2001]{Sun-Hu-Liu} by Zhi-Wei Sun, Hong Hu, J.-X. Liu  and \cite[2001]{Hu-Sun}  by Hong Hu and Zhi-Wei Sun). And with respect to natation:  If in  \cite[1878]{EdL} Fran\c{c}ois  \'Edouard Anatole Lucas  had been used $a=\textbf{p}$ and $b = \textbf{q}$ notation,  he would be perhaps at first glance notified and recognized as a Great Grandfather of all the $(p,q)$ - people. Let us start then introducing reconciling and matched denotations and nomenclature.

\begin{equation}
H_{n+2} = P \cdot H_{n+1}  - Q \cdot H_n ,\ \;n\geq 0 \ and \;  H_0 = a,\; H_1 = b.                                                                                     
\end{equation}
which is sometimes being written  in $\left\langle P,-Q \right\rangle \mapsto  \left\langle s,t \right\rangle$ notation.

\begin{equation}
H_{n+2} = s \cdot H_{n+1}  + t\cdot H_n ,\ \;n\geq 0 \ and \;  H_0 = a,\; H_1 = b.                                                                                     
\end{equation}

\vspace{0.2cm}
\noindent Simultaneously and collaterally  we  mnemonically pre adjust the starting point to discuss the $F(x)$ polynomials' case via - if entitled - antecedent "$\mapsto$ action": 
$H \mapsto H(x)$, $s \mapsto s(x)$, $t \mapsto t(x)$, etc. 

\begin{equation}
H_{n+2}(x) = s(x) \cdot H_{n+1}(x) + t(x)\cdot H_n, \; n\geq 0, \ H_0 = a(x), H_1 = b(x).                                                                                     
\end{equation}

\vspace{0.2cm}
\noindent enabling recovering  explicit formulas also for sequences of  polynomials correspondingly generated by the above linear recurrence of order 2  - with  Tchebysheff  polynomials and the generalized Gegenbauer-Humbert polynomials included. See for example Proposition 2.7 in the  recent Ward-Horadam peoples' paper  \cite [2009]{He-Siue}  by  Tian-Xiao He and Peter Jau-Shyong Shiue.

\vspace{0.2cm}
\noindent The general solution of (1):   $ H(a,b;P,Q) = \left\langle H_n\right\rangle _{n\geq 0}$ is being called  throught  this paper - \textbf{Ward-Horadam number'sequence}.

\vspace{0.3cm}

\noindent The general solution of (3):   $ H(x) \equiv H(a,b(x);s(x),t(x)) = \left\langle H_n(x)\right\rangle _{n\geq 0}$ is being called  throughout this paper - \textbf{Ward-Horadam functions' sequence}. It is then to be noted here that ideas germane to \textit{special  Ward-Horadam   polynomials  sequences}  of the \cite{1 Horadam 1996}  paper  were already explored in some details in  \cite{Horadam1}.  For more on special  Ward-Horadam   \textit{polynomials  sequences}  by Alwyn F. Horadam - consult  then:  \cite{Horadam 1969},  \cite[1985] {Horadam3},   \cite{Horadam 1993} ,  \cite{2 Horadam 1996} or  see \textit{for example} the following papers and references therein:   recent papers  \cite[2009]{HorzumKocer} by Tugba Horzum and  Emine  G{\"{o}}kcen Kocer  and  \cite [2009]{Gi-Sang Cheon}  by Gi-Sang Cheon, Hana Kim   and   Louis W. Shapiro.  For   \textit{Ward-Horadam   functions  sequences}    \cite[2009]{He-Siue}  by Tiang-Xiao He and Peter J. -S. Shiue who however there then concetrate on  on \textit{special  Ward-Horadam   polynomials  sequences} only.  

\vspace{0.1cm}

\noindent In \cite[2010]{Cigler q-Lucas 2010} Johann Cigler considers special Ward-Horadam   \textit{polynomials  sequences} and  among others he supplies the tiling combinatorial interpretation of these special Ward-Horadam   \textit{polynomials  sequences} which are  $q$-analogues of the Fibonacci and Lucas polynomials introduced in \cite[2002]{Cigler q-Fib q-Luc 2002}  and  \cite[2003]{Cigler Fib Pol 2003} by  Johann Cigler.

\vspace{0.1cm}

\noindent In the paper  \cite[2003]{Cigler} Johann Cigler  introduces "abstract Fibonacci polynomials"  - interpreted in terms of Morse coding sequences monoid with concatenation (monominos and dominos tiling then)   Cigler's abstract Fibonacci polynomial sare monoid algebra over reals valued polynomials with straightforward Morse sequences i.e. tiling recurrence originated (1.6)  "addition formula"

$$F_{m+n}(a,b) =  F_{m+1}(a,b)\cdot F_{m}(a,b) + b\cdot F_{n-1}(a,b)\cdot F_{n}(a,b),$$
which is attractive and seductive to deal with   within the context of this paper Theorem 1 below.                                                                                                                                                                                                                       

\vspace{0.2cm}

\noindent  From the characteristic equation of (1)

\begin{equation}
x^2 = P \cdot x  - Q,                                                                                      
\end{equation}
written by some of $p.q$-people   as

\begin{equation}
x^2 = s \cdot x  + t                                                                                      
\end{equation}
we readily find the Binet form solution of (1) (see (6) in \cite[2009]{HorzumKocer})  which is given by (6) and (7).

\begin{equation}
H_n(a,b;P,Q) \equiv H_n(A,B;p,q) = A p^n  + B q^n, \;n\geq 0 , H_0 = A, H_1 = B.                                                                                     
\end{equation}
where  $p,q$ are roots of (3) and we have assumed since now on that  $p\neq q$. 
As for the case $p=q$  included see for example Proposition 2.1 in  \cite[2009]{He-Siue} and see references therein.

\vspace{0.2cm}

\noindent Naturally :  $p+q = P\equiv s$  , $p \cdot q = Q \equiv -t$  and 

\begin{equation}
A = \frac{b-qa}{p-q} \ , \;    B = - \frac{b-pa}{p-q}.   
\end{equation}
\noindent hence we may and we shall use the following  conventional  identifications-abbreviations :

\begin{equation}
H \equiv H(a,b;P,Q) \equiv  H(A,B;P,Q) \equiv  H(A,B;p,q).
\end{equation}

\vspace{0.2cm}

\noindent It is obvious that the exponential generating function for Ward-Horadam sequence $H$ reads:

\begin{equation}
E_H(A,B;p,q)[x] = A exp[p\cdot x] + B exp[q\cdot x].   
\end{equation}

\vspace{0.2cm}

\noindent The derivation of the formula for ordinary generating function for Ward-Horadam  sequence is a standard task and so we have (compare with  (5) in \cite{HorzumKocer}) 

\begin{equation}
G_H(a,b;P,Q)[x] = \frac{a+(b-aP)x}{1-P \cdot x + Q \cdot x^2}=  \frac{a+(b-a[p+q])x}{1-P \cdot x + p\cdot q \cdot x^2}.   
\end{equation}
where from we decide an identification-abbreviation $$ G_H[x] \equiv  G_H(A,B;p,q)[x] \equiv G_H(a,b;P,Q)[x].$$

\vspace{0.2cm}

\noindent Naturally - in general  $H(A,B;p,q) \neq H(A,B;q,p)$. If $H(A,B;p,q) = H(A,B;q,p)$ we then call the  Ward-Horadam sequence symmetric and thus we arrive to Lucas
{\it Th\'eorie des Fonctions Num\'eriques Simplement Périodiques}  \cite[1878]{EdL}.

\vspace{0.2cm}

\noindent  In \cite[1878]{EdL} Edouard Lucas considers  Lucas sequence of the second kind   $V = \left\langle V_n\right\rangle_{n\geq 0}$ (second kind - see: \cite[1977, Part I]{KKK}) as well  as its till now well elaborated companion Lucas sequence of the first kind $ U = \left\langle U_n\right\rangle_{n\geq 0}$ (first  kind - see: \cite[1977, Part I]{KKK}) which gives rise in its turn to the  \textit{U}-binomials' recurrence (58) in \cite[1878]{EdL} (see then  \cite[1949]{JM},  \cite[1964]{TF},  \cite[1969]{Gould},  \cite[1989]{G-V}  or in  \cite[1989]{K-W}  etc. )

\vspace{0.2cm}

\noindent These sequences i.e ($A=B = 1$) the Lucas sequence of the second kind

\begin{equation}
H_n(2,P;p,q)= V_n = p^n + q^n.   
\end{equation}
\noindent and ($A= - B = 1$) the Lucas sequence of the first kind

\begin{equation}
H_n(0,1;p,q)= U_n = \frac{p^n - q^n}{p-q},    
\end{equation}
where called by Lucas \cite[1878]{EdL} the \textit{simply periodic numerical functions} because of 

\vspace{0.1cm}

\noindent [quote] \textsl{at the start, the complete analogy of these symmetric functions with 
the circular and hyperbolic functions.} [end of quote].

\vspace{0.3cm}


\noindent \textbf{More ad Notation 1.} 
\noindent The letters  a,b $a\neq b$ in \cite[1878]{EdL} denote the roots of the equation $x^2= P x -Q$  then  $(a,b)\mapsto (u,v)$ in \cite[2009]{Savage} and u,v stay there for the roots of the equation $x^2  =  \ell x -1$.

\noindent We shall use here the identification  $(a,b) \equiv (p,q)$  i.e.  $p,q$ denote the roots of $x^2= P x -Q$ as is common in "`Lucas $(p,q)$-people"' publications.


\noindent For Lucas $(p,q)$-people  then the following  $U$-identifications  are expediency natural:

\begin{defn}
     \begin{equation}
n_{p,q} = \sum_{j=0}^{n-1}{p^{n-j-1}q^j}=U_n =\frac{p^n - q^n} {p-q},\  0_{p,q}=U_0 = 0, \ 1_{p,q}=U_1 = 1,
    \end{equation}
\end{defn}
where  $p,q$  denote now the roots of the equation $x^2= P\cdot x - Q  \equiv  x^2= s x +t$ hence  $p+ q = s \equiv P$ , $pq = Q \equiv - t$ and the empty sum convention was used for $0_{p,q} = 0$.
\noindent Usually one assumes $p\neq q$. In general also $s\neq t$ - though  according to the context \cite[1989]{G-V}  $s=t$  may happen  to be the case of interest.

\vspace{0.2cm}

\noindent The Lucas $U$-binomial coefficients ${n \choose k}_U \equiv {n \choose k}_{p,q}$  are then defined as follows: (\cite[1878]{EdL},  \cite[1915]{Fon},  \cite[1936]{Ward},  \cite[1949]{JM},  \cite[1964]{TF}, \cite[1969]{Gould} etc.)

\vspace{0.3cm}

\begin{defn}
	Let  $U$ be as in \cite[1878]{EdL} i.e $U_n \equiv n_{p,q}$ then $U$-binomial coefficients for any $n,k \in \mathbb{N}\cup\{0\}$ are defined as follows
	\begin{equation}
		{n \choose k}_U \equiv {n \choose k}_{p,q} = \frac{n_{p,q}!}{k_{p,q}! \cdot (n-k)_{p,q}!} = \frac{n_{p,q}^{\underline{k}}}{k_{p,q}!}
	\end{equation}
	\noindent where $n_{p,q}! = n_{p,q}\cdot(n-1)_{p,q}\cdot ... \cdot 1_{p,q}$ and $n_{p,q}^{\underline{k}} = n_{p,q}\cdot(n-1)_{p,q}\cdot ...\cdot (n-k+1)_{p,q}.$
\end{defn} 

\vspace{0.3cm}

\begin{defn}
Let $V$ be as in \cite[1878]{EdL} i.e $V_n = p^n +q^n$, hence $V_0 = 2$ and  $V_n = p +q =s$. Then $V$-binomial coefficients for any $n,k \in \mathbb{N}\cup\{0\}$ are defined as follows
	\begin{equation}
		{n \choose k}_V =\frac{V_n!}{V_k!\cdot V_(n-k)!} = \frac{V_n^{\underline{k}}}{V_k!}
	\end{equation}
	\noindent where $V_n! = V_n \cdot V_{n-1}\cdot...\cdot V_1$ and $V_n^{\underline{k}}=V_n \cdot V_{n-1}\cdot ...\cdot V_{n-k+1}.$
\end{defn}

\vspace{0.2cm}

\noindent One automatically generalizes number $F$-binomial coefficients' array to functions $F(x)$-\textbf{multinomial} coefficients' array (see  \cite[2004]{akk3 2004} and references to umbral calculus therein) while for \textit{number sequences} $F$ the $F = F(x=1)$-multinomial coefficients see p. 40 in \cite[1913]{Carmichel} by  Robert Daniel Carmichel , see  \cite[1936]{Ward} by Morgan Ward,  \cite[1969]{Gould}  by Henri W. Gould or \cite[1979]{Shannon 1979 multi}) by  Anthony G. Shannon where recursion for $U = F(x=1)$-multinomial coefficients is provided and where there specifically $U$ denotes the fundamental Lucas sequence (i.e. the Lucas sequence of the first kind) and see also important paper \cite[1991]{SatoShi} by Shiro Ando and  Daihachiro Sato.  The  $x$-Fibonomial coefficients  from  \cite[2001]{Rich 2001}  by  Thomas M. Richardson are motivating example of  functions $F(x)$-\textbf{binomial} coefficients' array from \cite[2004]{akk3 2004}.

\begin{defn}\label{def:symbol}
Let $F(x)$ be any natural, or complex  numbers' non zero valued functions' sequence i.e. $F_n(x)\in\mathbb{N}$ or and $F_n(x)\in\mathbb{C}$. The \textbf{$F(x)$-multinomial coefficient} is then identified with the symbol

\begin{equation}
	{n \choose {k_1,k_2,...,k_s}}_{F(x)} = \frac{F_n(x)!}{F_{k_1}(x)!\cdot ... \cdot F_{k_s}(x)!}
\end{equation}

\vspace{0.2cm}
\noindent where $k_i\in\mathbb{N}$ and $\sum_{i=1}^{s}{k_i} = n$ for $i=1,2,...,s$. Otherwise it is equal to zero,
\end{defn}
\noindent and where $F_r(x)! = F_r(x) \cdot F_{r-1}(x)\cdot...\cdot F_1(x)$.

\vspace{0.2cm}

\noindent Naturally  for any natural $n,k$ and $k_1+...+k_m=n-k$ the following holds

\begin{equation} 
{n \choose k}_{F(x)} \cdot {n-k \choose {k_1,k_2,...,k_m}}_{F(x)}  = {n \choose {k,k_1,k_2,...,k_m}}_{F(x)}, 
\end{equation}
 
$${n \choose {k_1,k_2,...,k_m}}_{F(x)}={n \choose k_1}_{F(x)}{n-k_1\choose{k_2}}_{F(x)}\cdots {n-k_1- \cdots -k_{m-1}\choose{k_m}}_{F(x)}.$$


\vspace{0.3cm}


\noindent \textbf{More ad Notation 2.} 

\noindent We shall use further on the traditional , XIX-th century rooted notation under presentation in spite of being inclined to quite younger notation from \cite[2010]{BSCS}  by  Bruce E. Sagan and Carla  D.  Savage.  This wise, economic notation is ready for straightforward 	record of combinatorial interpretations and combinatorial interpretations' substantiation in terms of popular text book tiling model since long ago used for example to visualize recurrence for Fibonacci-like sequences ; see for example \cite[1989]{GKP 1989}  by Ronald Graham, Donald Ervin  Knuth, and Oren Patashnik.  
The translation from Fran\c{c}ois  \'Edouard Anatole Lucas via Dov Jarden and Theodor Motzkin  notation \cite[1949]{JM} and notation  of Bruce E. Sagan and Carla  D.  Savage \cite[2010]{BSCS} is based on the succeeding identifications: the symbol used for $U$-\textbf{binomials} is  $C\left\{... \right\}$  in place of  $\left(...\right)_U$ , the would be symbol  for $V$-\textbf{binomials} i.e. $P\left\langle ... \right\rangle$  in place of  $\left(...\right)_V$  is not considered at all in  \cite{BSCS} while  

$$ \left\{n \right\} \equiv U_n \equiv n_{p,q}  , \ \left\langle n \right\rangle \equiv V_n.$$
In Bruce E. Sagan and Carla  D.  Savage  notation  we would then write down the  fundamental and primordial sequences' binomial coefficients  as follows.

\vspace{0.3cm}

\begin{defn}
	Let  $\left\{n\right\}$ be fundamental Lucas sequence as in \cite[1878]{EdL} i.e $\left\{n\right\}\equiv U_n \equiv n_{p,q}$ then $\left\{n\right\}$-binomial coefficients for any $n,k \in \mathbb{N}\cup\{0\}$ are defined as follows
	\begin{equation}
		F\left\{n,k \right\} = \left\{ n \atop k\right\}_{p,q} = \frac{\left\{n\right\}!}{\left\{k\right\}! \cdot \left\{n-k\right\}!} = \frac{\left\{n\right\}^{\underline{k}}}{\left\{k\right\}!}
\end{equation}
where $\left\{n \right\}! = \left\{n \right\}\cdot\left\{n-1\right\}\cdot...\cdot \left\{1\right\}$  
and $\left\{n \right\}^{\underline{k}} = \left\{n\right\}\cdot\left\{n-1\right\}\cdot ...\cdot\left\{n-k+1\right\},$
\end{defn}

\begin{defn}
	Let  $\left\langle n \right\rangle $ be primordial Lucas sequence as in \cite[1878]{EdL} i.e $\left\langle n \right\rangle \equiv V_n $ then $ \left\langle n \right\rangle $-binomial coefficients for any $n,k \in \mathbb{N}\cup\{0\}$ are defined as follows
	\begin{equation}
P\left\langle n,k \right\rangle = \left\langle  n \atop k\right\rangle_{p,q} =\frac{\left\langle n \right\rangle!}{\left\langle k \right\rangle! \cdot \left\langle n-k \right\rangle!} = \frac{\left\langle n \right\rangle^{\underline{k}}}{\left\langle k \right\rangle!},
\end{equation}
where $\left\langle n \right\rangle! = \left\{n \right\}\cdot\left\{n-1\right\}\cdot...\cdot \left\{1\right\}$  
and $\left\{n \right\}^{\underline{k}} = \left\{n\right\}\cdot\left\{n-1\right\}\cdot ...\cdot\left\{n-k+1\right\}.$
\end{defn} 
The above consequent symbols $\left\{ n \atop k\right\}_{p,q}$  and  $\left\langle  n \atop k\right\rangle_{p,q}$     are - in  not exceptional  conflict - with second kind Stirling numbers notation and Euler numbers notation respectively in the spirit of \cite[1989]{GKP 1989} what extends on both  $p,q$ - extensions' notation.  

\noindent Regarding  the symbol $\left\{ n \atop k\right\}_{p,q}$ one draws the attention of a reader to  \cite[1967]{Hoggatt 1967} where Verner Emil Hoggatt, Jr.    considers the $C$-binomial coefficients with indices in an arithmetic progression  denoting them by symbols 
$ \left\{ n \atop k\right\}_{u,k}$, where $\left\{u_n\right\}_{n\geq 0}= U $  with  $U$being the primordial Lucas sequence.  For  $ \left\{ n \atop k\right\}$ corresponding notation see also: \cite[1969]{Gould}  by Henri W.Gould, \cite[1989]{G-V}   by Ira  M. Gessel and Xavier G\'erard Viennot , \cite[2005]{ST 2005}, \cite[2005]{ST 2005 2}, \cite[2006]{ST 2006}  by Jaroslav Seibert  and Pavel Trojovsk\'y and \cite[2007]{T 2007}  by  Pavel Trojovsk\'y.

\vspace{0.2cm}

\noindent Whereas as in the\textit{ subset-subspace problem} ( Example [Ex. q* ; 6] in subsection 4.3.)  we rather need another natural notation. Namely  for  $q \neq 0$  introduce  $q* = \frac{p}{q}$ and observe that

$$ {n \choose k}_U \equiv {n \choose k}_{p,q} = q^{k(n-k)} {n \choose k}_{1,q*} \stackrel{q*\mapsto 1}{\rightarrow}   \fnomialF{n}{k}. $$    
The $V$-binomial  $P\left\langle n,k \right\rangle = \left\langle  n \atop k\right\rangle_{p,q} \equiv \fnomialF{n}{k}{V}$ is not considered in \cite[2010]{BSCS}.




\section{$H(x)$-binomial coefficients' recurrence}

\vspace{0.2cm}

\noindent \textbf{3.1.}  Let us recall convention resulting from (3).

\vspace{0.2cm}
\noindent Recall. The general solution of (3):   $ H(x) \equiv H(a,b(x);s(x),t(x)) = \left\langle H_n(x)\right\rangle _{n\geq 0}$ is being called  throughout this paper - \textbf{Ward-Horadam functions' sequence}.

\vspace{0.2cm}

\noindent  From the characteristic equation of the recurrence (3)

\begin{equation}
z^2 - s(x) \cdot z  - t(x) =0
\end{equation}
we readily see that  for $ H_0 = a(x), \:H_1 = b(x), \;n\geq 0 ,$

\begin{equation}
H_n(x) \equiv H_n(a(x),b(x);p(x),q(x)) = A(x)p(x)^n + B(x)q(x)^n,                                                                                     
\end{equation}
where  $p(x),q(x)$ are roots of (20) and we have assumed that  $p(x)\neq q(x)$  as well as that  $p(x), q(x)$ are not roots of unity. Naturally:
 
\begin{equation}
A(x) = \frac{b(x)-q(x)a(x)}{p(x)-q(x)} \ , \;    B = - \frac{b(x)-p(x)a(x)}{p(x)-q(x)}.   
\end{equation}
\noindent hence we may and we shall use the following  conventional  identifications-abbreviations

\begin{equation}
H(x) \equiv H(a(x),b(x);s(x),t(x)) \equiv  H(A(x),B(x);s(x),t(x)).
\end{equation}

\vspace{0.2cm}

\noindent As for the case $p(x)=q(x)$  included see for example Proposition 2.7 in  \cite[2009]{He-Siue}.

\vspace{0.2cm}

\noindent Another  explicit formula for Ward-Horadam functions sequences is the mnemonically  extended formula (9) from  \cite[2009]{HorzumKocer} by Tugba Horzum and Emine G{\"{o}}kcen  Kocer,  where here down we use contractually the following abbreviations:

\vspace{0.2cm}

\noindent$H_n(x) \equiv H_n(a(x),b(x);s(x),t(x))$ , $a(x) \equiv a $, $b(x) \equiv b $, $s(x) \equiv s $ and  $t(x) \equiv t $

\begin{equation}
H_n(x)= a\sum_{0\leq k \leq \left\lfloor \frac{n}{2}\right\rfloor} \binom{n-k}{k}s^{n-2k}t^k + \left( \frac{b}{s}-a \right)\sum_{0\leq k \leq \left\lfloor \frac{n-1}{2}\right\rfloor} \binom{n-k-1}{k}s^{n-2k}t^k                                                                          
\end{equation}

\vspace{0.2cm}

\noindent \textbf{Note and compare}. The recurrence  (1.1) and (1.2) in \cite[1996]{1 Horadam 1996} defines a polynomials' subclass of Ward-Horadam functions sequences defined by (3).
The standard Jacques Binet  form (1.8) in \cite[1996]{1 Horadam 1996} of the  recurrence  (1.1) and (1.2) solution for Ward-Horadam polynomials  sequences in \cite[1996]{1 Horadam 1996} is the standard Jacques Binet  form (19), (20) of the  recurrence  (3) solution for Ward-Horadam functions sequences.

\vspace{0.1cm}

\noindent The recurrence  (2.23) in \cite [2009]{He-Siue}  by  Tian-Xiao He and   Peter Jau-Shyong Shiue defines exactly the class of Ward-Horadam functions' second order sequences 
and the standard Jacques Binet  form (19), (20) of the  recurrence  (3) solution for Ward-Horadam functions sequences $H(x)$ constitutes the content of their Proposition 2.7. - as has been  mentioned earlier.  No recurrences for $H(x)$-binomials neither for $H(x=1)$-binomials are considered.

\vspace{0.2cm}

\noindent \textbf{On Binet Formula - Historical Remark.}  We just quote Radoslav Rasko Jovanovic's information from   $$ http://milan.milanovic.org/math/english/relations/relation1.html :$$

\begin{quot} Binet's Fibonacci Number Formula was derived by Binet in 1843 although the result was known to Euler and to Daniel Bernoulli more than a century ago. ... It is interesting that A de Moivre (1667-1754) had written about Binet`s Formula, in 1730, and had indeed found a method for finding formula for any general series of numbers formed in a similar way to the Fibonacci series.
\end{quot} 
See also the book  \cite[1989]{Vajda 1989}    by  Steven Vajda. 

\vspace{0.2cm}

\noindent \textbf{3.2.}  The authors of \cite{Savage} provide an easy proof of  an observation named there Theorem 17  which extends automatically to the statement that  the following  recurrence holds for the general case of ${r+s \choose r,s}_{H(x)}$  $H(x)$-binomial arrays in multinomial notation.

\vspace{0.2cm}

\begin{theoremn}
Let us admit shortly the abbreviations: $g_k(r,s)(x) = g_k(r,s)$ , $k=1,2$.  Let $s,r>0$. Let $F(x)$ be any zero characteristic field nonzero valued functions' sequence ($F_n(x) \neq 0$). Then

\begin{equation}
{r+s \choose r,s}_{F(x)}=g_1(r,s)\cdot{r+s-1 \choose r-1,s}_{F(x)}+g_2(r,s)\cdot{r+s-1 \choose r,s-1}_{F(x)}
\end{equation}
where   $ {r \choose r,0}_{F(x)} = {s \choose 0,s}_{F(x)} =1$ and 
\begin{equation}
 F(x)_{r+s} =  g_1(r,s) \cdot F(x)_r   +  g_2(r,s) \cdot F(x)_s.                               
\end{equation}
are equivalent.
\end{theoremn}

\vspace{0.3cm}

\noindent \textbf{On the way historical note}
\noindent  Donald  Ervin  Knuth   and    Herbert Saul Wilf in \cite[1989]{K-W} stated that  Fibonomial coefficients and the recurrent relations for them appeared already in 1878  Lucas work (see: formula (58)  in \cite[1878]{EdL}  p. 27 ;  for $U$-binomials which "Fibonomials" are special case of). More over on this very p. 27 Lucas formulated a conclusion from his (58) formula which may be  stated in notation of this paper formula (2) as follows: \textit{if}  $s,t  \in\mathbb{Z}$ and $H_0 =0$ , $H_1 = 1$  \textit{then} $H \equiv U$  and $ \fnomialF{n}{k}{U} \equiv \fnomialF{n}{k}{n_{p,q}} \in\mathbb{Z}$.  Consult also in next century  \cite[1910]{Bachmann 1910}  by  Paul Gustav Heinrich Bachmann  or later on - \cite[1913]{Carmichel}  by  Robert Daniel Carmichel  [p. 40]  or  \cite[1949]{JM} by Dov Jarden and Theodor Motzkin where in all quoted positions it was also shown that $n_{p,q}$ - binomial coefficients are integers  - for $p$ and $q$ representing  distinct roots of (5) with  their ratio being not a root of unity.

\noindent Let us take an advantage to note that  Lucas  Th\'eorie des Fonctions Num\'eriques Simplement P\'eriodiques i.e. investigation exactly of fundamental $U$  and  primordial $V$ sequences constitutes  the far more non-accidental context  for binomial-type coefficients exhibiting their relevance at the same time to number theory and  to hyperbolic trigonometry  (in addition to \cite[1878]{EdL} see for example \cite{akk1999},  \cite{ Bakk} and  \cite{Bajguz}).  

\vspace{0.1cm}

\noindent It seems to be the right place now to underline that the  \textit{addition formulas} for Lucas sequences below with respective hyperbolic trigonometry formulas and also consequently $U$-binomials'recurrence formulas - stem from commutative  ring $R$ identity: $(x-y)\cdot (x+y) \equiv x^2 - y^2, x,y \in R$.

\vspace{0.2cm}

\noindent Indeed. Taking here into account the  $\textbf{U}$\textit{-addition formula} i.e. the first of two trigonometric-like $L$-addition formulas (42) from \cite[1878]{EdL} ($L[p,q]  = L = U,V$ - see also \cite[1999]{akk1999} by  A.K.Kwa\'sniewski and   \cite{ Bakk}, \cite{Bajguz})  i.e.

\begin{equation}
2 U_{r+s} =  U_r V_s  +  U_s V_r ,\ \ \ \ 
2 V_{r+s} =  V_r V_s  +  U_s U_r 
\end{equation}
one readily recognizes that  the $U$-binomial recurrence from the Corollary 18 in  \cite[2009]{Savage} is a case of the $U$-binomial recurrence (58) \cite[1878]{EdL}  which may be rewritten after Fran\c{c}ois  \'Edouard Anatole Lucas in multinomial notation and stated as follows: \textit{according  to the Theorem 1  case i.e. the Theorem 2 below the following is true}

$$2 U_{r+s} =  U_r V_s  +  U_s V_r$$ 
\textit{is equivalent to }

\begin{equation}
2\cdot{r+s \choose r,s}_{n_{p,q}}= V_s \cdot{r+s-1 \choose r-1,s}_{n_{p,q}}+V_r \cdot{r+s-1 \choose r,s-1}_{n_{p,q}}.
\end{equation}
\noindent  To this end see also Proposition 2.2. in \cite[2010]{BSCS} and compare it with both (28) and  Example 3. below.

\vspace{0.2cm}

\noindent However there is no companion  $V$-binomial recurrence i.e. for ${r+s \choose r,s}_V$ neither in  \cite[1878]{EdL}  nor in \cite[2009]{Savage} as well as all other quoted papers except for \cite[2010]{akk2010}  - up to knowledge of this note author.  

\vspace{0.2cm}

\noindent Consequently  then there is no  $H(x)$-binomial recurrence neither in  \cite[1878]{EdL}  nor in \cite{Savage} (2009) as well as all other quoted papers  except for Final remark : p.5 in \cite[2010]{akk2010} up to this note author knowledge.  

\vspace{0.2cm}

\noindent \textbf{The End} of \textit{the on the way historical note}.

\vspace{0.3cm}

\noindent The looked for  $H(x)$-binomial recurrence  (29) accompanied by (30-33) might be  then given right now in the form of (25) adapted to  - Ward-Lucas functions'sequence case notation while keeping in mind that of course the expressions for $h_k(r,s)(x)$,  $k=1,2$ below are designated by this $ F(x) = H(x)$ choice and as a matter of fact are appointed by the recurrence (3). 

\vspace{0.3cm}

\noindent For the sake of commodity let us admit shortly the abbreviations: $h_k(r,s)(x) = h_k(r,s)= h_k$ , $k=1,2$. Then for $H(x)$ of the form (21) we evidently have what follows.

\vspace{0.4cm}

\noindent \textbf{Theorem 2.}

\begin{equation}
{r+s \choose r,s}_{H(x)} = h_1(r,s){r+s-1 \choose r-1,s}_{H(x)} + h_2(r,s){r+s-1 \choose r,s-1}_{H(x)},
\end{equation}
where $p(x) \neq q(x)$ and  $ {r \choose r,0}_{H(x)} = {s \choose 0,s}_{H(x)} =1,$  \textit{is equivalent to}

\begin{equation}
 H(x)_{r+s} =  h_1(r,s)H(x)_r   +  h_2(r,s) H(x)_s.                               
\end{equation}
where  $H_n(x)$ is explicitly given by (21)  and (22). 
\vspace{0.1cm}
\noindent \textbf{The end} of the Theorem 2.

\vspace{0.2cm}

\noindent There might be various  $h_1(r,s)(x)=h_1$ and $ h_2(r,s)(x)_s =h_2 $ solutions of  (30) and  (21). Compare (38)  in Example 1 with (42) in Example 3  below. 
\noindent As the possible  $h_1(r,s)(x)=h_1$ and $ h_2(r,s)(x)_s =h_2 $ formal solutions of  (30) and  (21)  we may take

\begin{equation}
h_1(r,s)(x)= \frac{A(x)\cdot p(x)^{r+s}} {A(x) \cdot p^r + B(x) \cdot q(x)^r},   
\qquad
h_2(r,r)(x) = \frac{B(x) \cdot q(x)^{r+s}} {A(x) \cdot p^s + B(x) \cdot q(x)^s}.
\end{equation}
As  another possible $h_1(r,s)(x)$ and $ h_2(r,s)(x)_s =h_2 $ solutions of  (30) and  (21) we may take:   \textbf{for $r \neq s$}

\begin{equation}
h_1(r,s) \cdot (p((x)^r q(x)^s  -  q(x)^r p(x)^s ) = p(x)^{r+s} q(x)^s  -  q(x)^{r+s} p(x)^s,
\end{equation}

\begin{equation}
h_2(r,s) \cdot (q(x)^r p(x)^s  -  p(x)^r q(x)^s ) = p(x)^{r+s} q(x)^r  -  q(x)^{r+s} p(x)^r. 
\end{equation}
\textit{while for }$r=s$  apply formula (31)  with $r=s$.

\vspace{0.2cm}

\noindent Usually the specific features of  particular cases of (21)  and (22) allow one to infer the particular form of (30) hence the form of $h_1(r,s)(x)=h_1$ and $ h_2(r,s)(x)_s =h_2 $.

\vspace{0.5cm}

\noindent \textbf{3.3. Three special cases examples.}

\vspace{0.2cm}

\noindent \textbf{Example 1.} This is a particular case of the Theorem 2. 
\vspace{0.1cm}

\noindent The recurrent relations (13) and (14) in  Theorem 1 from  \cite[2008]{Corsino} by Roberto Bagsarsa Corcino  for $n_{p,q}$-binomial coefficients are special cases of this paper formula (29) as well as  of Th. 17 in \cite{Savage} with straightforward identifications of $g_1, g_2$  in (13)  and  in (14) in \cite{Corsino} or in this paper recurrence (30) for $H(x=1) = U[p,q]_n = n_{p,q}$ sequence. Namely, recall here now in multinomial notation this  Theorem 1 from  \cite[2008]{Corsino}   by Roberto Bagsarsa Corcino:

\begin{equation}
{r+s \choose r,s}_{p,q} = q^r{r+s-1 \choose r-1,s}_{p,q} + p^s{r+s-1 \choose r,s-1}_{p,q},
\end{equation}

\begin{equation}
{r+s \choose r,s}_{p,q} = p^r{r+s-1 \choose r-1,s}_{p,q} + q^s{r+s-1 \choose r,s-1}_{p,q},
\end{equation}
which is equivalent to

\begin{equation}
(s + r)_{p,q} = p^s r_{p,q} +  q^r s_{p,q} = (r + s)_{q,p} = p^r s_{p,q} +  q^s  r_{p,q},
\end{equation}
what  might be at once seen proved by noticing that $$p^{r+s} - q^{r+s} \equiv p^s \cdot(p^r - q^r) +   q^r \cdot(p^s - q^s).$$  Hence those mentioned straightforward identifications follow:

\begin{equation}
 g_1 = q^r ,\ \  g_2 =  p^s \;or\ \  g_1 = p^r ,\ \  g_2 =  q^s. 
\end{equation}
The recurrence (36) in Lucas notation reads
\begin{equation}
U_{s+r} = p^s U_r +  q^r U_s = U_{r+s} = p^r U_s +  p^s U_r. 
\end{equation}
Compare it with equivalent recurrence (42) in order to notice that  both $h_1$  and $h_2$ functions are different from case to case  of recurrence (30) \textbf{equivalent realizations}.

\vspace{0.2cm}

\noindent  Compare this example based on Theorem 1  in \cite[2008]{Corsino} by Roberto Corcino with  with \cite[2008]{MD2} \textbf{v[1]} by Maciej Dziemia\'nczuk (see there (1) and (2) formulas), and track as well  -  the simple combinatorial proof of the Corollary 3   in  \cite[2009]{MD2} \textbf{v[2]})  by Maciej Dziemia\'nczuk.

\vspace{0.3cm}

\noindent \textbf{Example 2.} This is a particular case of the Theorem 1. 
\vspace{0.1cm}

\noindent Now let $A$ be any natural numbers' or even complex numbers' valued sequence. One readily sees that also (1915) Fonten\'e recurrence for Fonten\'e-Ward generalized $A$-binomial coefficients i.e. equivalent identities (6) , (7) in \cite{Gould} \textbf{are special cases of} this paper formula \textbf{(26)} as well as  of Th. 17 in \cite{Savage} with straightforward identifications of $h_1, h_2$  in this paper formula (25)  while this paper recurrence  (27) becomes trivial identity.

\noindent Namely, the identities (6) and (7) from \cite[1969]{Gould} read correspondingly:

\begin{equation}
{r+s \choose r,s}_A = 1 \cdot {r+s-1 \choose r-1,s}_A + \frac{A_{r+s} - A_r}{A_s}{r+s-1 \choose r,s-1}_A,
\end{equation}

\begin{equation}
{r+s \choose r,s}_A = \frac{A_{r+s} - A_s} {A_{r}} \cdot {r+s-1 \choose r-1,s}_A +  1 \cdot {r+s-1 \choose r,s-1}_A,
\end{equation}
where $p \neq q$ and  $ {r \choose r,0}_L = {s \choose 0,s}_L =1.$ And finally we have tautology identity

\begin{equation}
A_{s + r} \equiv \frac{A_{r+s} - A_s}{A_r}\cdot A_r +  1 \cdot A_s.
\end{equation}
Example 2. becomes the general case of the Theorem 1. if we allow $A$ to represent any zero characteristic field nonzero valued functions' sequence: $A = A(x) = \left\langle A_n(x)\right\rangle_{n \geq 0}, \;  A_n(x) \neq 0$).

\vspace{0.3cm}

\noindent \textbf{Example 3.} This is a particular case of the Theorem 2. 
\vspace{0.1cm}

\noindent The first example above is cognate to this third example in apparent way  as might readily seen from Fran\c{c}ois  \'Edouard Anatole Lucas papers \cite[1878]{EdL} or more  recent article  \cite[2001]{Hu-Sun}  by Hong Hu and Zhi-Wei Sun  ; (see also $t=s$ case in \cite[1989]{G-V} by Ira  M. Gessel and Xavier G\'erard Viennot on pp.23,24 .) In order to experience  this let us start to consider now the number $H(x=1)=U$ Lucas fundamental sequence \textit{fulfilling} (\textbf{2}) with  $U_0 = 0$  and $U_1 = 1$  as  introduced in \cite[1878]{EdL} and the - for example  considered in \cite[2001]{Hu-Sun}. There in \cite[2001]{Hu-Sun} by Hong-Hu and Shi-Wei Sun - as a matter of fact - a kind of "pre-Theorem 17"  from \cite[2009]{Savage} is latent in the  proof of  Lemma 1 in  \cite{Hu-Sun}.  We rewrite  Lemma 1  by Hong-Hu and Shi-Wei Sun  in multinomial notation and an arrangement convenient for our purpose here using sometimes abbreviation $U_n(p,q)\equiv U_n$. 

\noindent (Note that the  \textit{addition formulas} for Lucas sequences hence consequently $U$-binomials'recurrence formulas  \cite[1878]{EdL} as well as  $(p-q)\cdot(p^{j+k} - q^{j+k}) \equiv (p^{k+1} - q^{k+1})\cdot(p^j - q^j) - p\cdot q (p^{j-1} - q^{j-1} \cdot(p^k - q^k)$ - stem from commutative ring $R$ identity: $(x-y)\cdot (x+y) \equiv x^2 - y^2, x,y \in R$.)  

\noindent And so for   $p \neq q$ and bearing in mind that  $p\cdot q = - t$ - the following is true. 

\vspace{0.2cm}

\noindent The identity (42) equivalent to 
$$(p-q)\cdot(p^{j+k} - q^{j+k}) \equiv (p^{k+1} - q^{k+1})\cdot(p^j - q^j) - p\cdot q (p^{j-1} - q^{j-1} \cdot(p^k - q^k)$$  

\begin{equation}
U_{j + k}(p,q) = U_{k+1} \cdot U_j(p,q) + t  U_{j-1}\cdot U_k(p,q)
\end{equation}
\noindent\textit{is equivalent to}

\begin{equation}
{j+k \choose j,k}_U = U_{k+1} \cdot {j+k-1 \choose j-1,k}_U + U_{j-1}\cdot{j+k-1 \choose j,k-1}_U,
\end{equation}
\textit{where} $p,q$  \textit{are the roots of} (5) \textit{and correspondingly the above Lucas fundamental sequence} $H_n = U_n(p,q)$ i.e.  $U_0 = 0$  and $U_1 = 1$ \textit{is given by its Binet form} (6),(7). 

\vspace{0.2cm}

\noindent Compare (42) with equivalent recurrence (48) in order to notice that  both $h_1$  and $h_2$ functions are different from case to case  of recurrence (30) \textbf{equivalent realizations}.  

\vspace{0.2cm}

\noindent Compare now: this paper recurrence formula (42)  with recurrence formula (4) in \cite[2010]{BSCS}, compare this paper recurrence formula (43)  with Proposition 2.2.  in \cite[2010]{BSCS}  by  Bruce E. Sagan and  Carla  D.  Savage. Compare this paper recurrence (28) equivalent to (5) and proposition 2.2. in \cite[2010]{BSCS}  and note that (5) in \cite[2010]{BSCS} is just the same - as (58) in \cite[1878]{EdL} - the same except for notation.  The translation from "younger" notation of Bruce E. Sagan and Carla  D.  Savage (from one - left hand - side) into more matured by tradition  notation of Fran\c{c}ois  \'Edouard Anatole Lucas (from the other - right hand - side) is based on the identifications: the symbol used for $U$-binomials is  $\left\{... \right\}$  in place of  $\left(...\right)_U$  and 

$$ \left\{ n \right\} \equiv U_n \equiv n_{p,q}  , \ \left\langle n \right\rangle \equiv V_n.$$

\vspace{0.2cm}

\noindent For $s=t=1$  we get Fibonacci $U_n = F_n$ sequence with recurrence (41) becoming the recurrence known from Donald Ervin Knuth and Herbert Saul Wilf  masterpiece \cite[1989]{K-W}.

Example 3. becomes more general case of the Theorem 1. if we allow $U$ to represent any zero characteristic field nonzero valued functions' sequence: $U(x) = \left\langle U_n(x)\right\rangle \_{n\geq 0}, \;  U_n(x) = \frac{p(x)^n-q(x)^n}{p(x)-q(x)} \equiv n_{p(x),q(x)},  p(x) \neq  q(x)$, where  $p(x),q(x)$ denote the distinct roots  of (20) and we have assumed  as well  that  $p(x), q(x)$ are not roots of unity. 

\vspace{0.2cm}

\noindent \textbf{The End} of three examples.

\vspace{0.2cm}


\section{Snatchy information on $F$-binomials' and their relatives' combinatorial interpretations }

\vspace{0.1cm}

\noindent \textbf{ 4.1.} In regard to  \textbf{combinatorial interpretations} of $L$-binomial or $F$-multinomial coefficients or related arrays we leave that subject apart from this note.  Nevertheless we direct the reader to some comprise papers and references therein; these are \textsl{for example} here the following:

\vspace{0.2cm}

\noindent Listing. \textbf{1.} \cite[1984]{Voigt} by Bernd Voigt: on common generalization of binomial coefficients, Stirling numbers and Gaussian coefficients .

\vspace{0.2cm}

\noindent Listing. \textbf{2.}  \cite[1991]{wachs} by Michelle L. Wachs and   Dennis White and in \cite[1994]{wachs 2} by Michelle L. Wachs: on p,q-Stirling numbers and set partitions.

\vspace{0.2cm}

\noindent Listing. \textbf{3.} \cite[1993]{medicis} by Anne De M\'edicis  and Pierre Leroux:  on Generalized Stirling Numbers, Convolution Formulae and (p,q)-Analogues.

\vspace{0.2cm}

\noindent Listing. \textbf{4.} \cite[1998]{Konva 1998} John Konvalina:  on generalized binomial coefficients and the  Subset-Subspace Problem.  Consult examples [Ex. q* ; 6] and  [Ex. q* ; 7] in \textbf{4.3}. below.  Then see also \cite[2000]{Konva 2000}  by   John Konvalina on an unified  simultaneous  interpretation of binomial coefficients of both kinds, Stirling numbers of both kinds and Gaussian binomial coefficients  of both kinds.

\vspace{0.2cm}

\noindent Listing. \textbf{5.} Ira  M. Gessel  and Xavier G\'erard Viennot in  \cite[1989]{G-V} deliver now well known their interpretation of the fibonomials in terms of non-intersecting lattice paths .

\vspace{0.2cm}

\noindent Listing. \textbf{6.} In \cite[2004]{RW}  Jeffrey B. Remmel and Michelle L. Wachs derive a new rook theory interpretation of a certain class of generalized
Stirling numbers and their $(p,q)$-analogues.  In particular they prove that their $(p,q)$-analogues of the generalized Stirling numbers of the second
kind may be interpreted in terms of colored set partitions and colored restricted growth
functions.

\vspace{0.2cm}

\noindent Listing. \textbf{7.}  \cite[2005]{Otta vio Munari}  by Ottavio M. D'Antona  and Emanuele Munarini deals with  - in terms of weighted binary paths - combinatorial interpretation of the connection constants which is in particular unified, simultaneous combinatorial interpretation  for Gaussian coefficients, Lagrange sum, Lah numbers, ,  q-Lah numbers, Stirling numbers of both kinds , q-Stirling numbers of both kinds. Notr the correspondence: weighted binary paths $\Leftrightarrow$     edge  colored binary paths

\vspace{0.2cm}
\noindent  Listing. \textbf{8.} Maciej  Dziemia\'nczuk in \cite[2011]{MD 2011}  extends the results of John Konvalina from \textbf{4.} above. The Dziemia\'nczuk'  $\zeta$ - analogues of the Stirling numbers arrays of both kinds cover ordinary binomial and Gaussian coefficients, $p,q$-Stirling numbers and other combinatorial numbers studied with the help of object selection, Ferrers diagrams and rook theory. The $p,q$-binomial arrays are special cases of $\zeta$- numbers' arrays, too. 

\vspace{0.1cm}
\noindent $\zeta$ -number of the  first and the second kind is the number of ways to select $k$  objects from $k$  of $ n$ boxes without box repetition allowed and with box repetition allowed, respectively.  

\vspace{0.1cm}
\noindent The weight vectors used for objects constructions and  statements derivation are functions of parameter $\zeta$.

\vspace{0.2cm}
\noindent  Listing. \textbf{9.} As regards combinatorial interpretations via  tilings in   \cite[2003]{BQ2003} and  \cite[2010]{BP}  - see \textbf{4.2.} below.

\vspace{0.2cm}

\noindent Listing. \textbf{10.} In   \cite[2003]{Cigler} Johann Cigler  introduces "abstract Fibonacci polynomials"  - interpreted in terms of Morse coding sequences monoid with concatenation (monominos and dominos tiling then).  Cigler's abstract Fibonacci polynomial sare monoid algebra over reals valued polynomials with straightforward Morse sequences i.e. tiling recurrence originated (1.6)  "addition formula"

$$F_{m+n}(a,b) =  F_{m+1}(a,b)\cdot F_{m}(a,b) + b\cdot F_{n-1}(a,b)\cdot F_{n}(a,b),$$
which is attractive and seductive to deal with   within the context of this paper Theorem 1. The combinatorial tiling interpretation of the model is its construction framed in the Morse coding sequences monoid with concatenation (monominos and dominos tiling then).

\vspace{0.2cm}

\noindent  Listing. \textbf{11.}  In \cite[2010]{Cigler q-Lucas 2010} Johann Cigler considers special Ward-Horadam   \textit{polynomials  sequences} and reveals the tiling combinatorial interpretation of these special Ward-Horadam   \textit{polynomials  sequences} in the spirit of Morse with monomino, domino alphabet monoid as here above in \textbf{10.}. Namely:

\vspace{0.1cm}

\noindent \textbf{1.} the $q$-Fibonacci polynomial $F_n(x,s,q) = \sum_{c\in \Phi_n}w(c)\equiv w(\Phi_n)$ is the weight function of the set $\Phi_n$ of all words  (coverings) $c$ of length $n-1$ in Morse (tiling) alphabet $\left\{a,b\right\}$ i.e. corresponding generation function for number of linear tilings as  $\Phi_n$ clearly with the set of may be identified with the set of all linear tilings  of  $(n-1)\times 1$  rectangle or equivalently with Morse code sequences of length $n-1$. 

\noindent Polynomials $F_n(x,s,q)$ satify this paper recursion (3)  with $H_0(x)=0$ , $H_1(x)=1$ ; $s(x)=x$ and $t(x)=s$.

\noindent The $F_n(x,s,q)$-binomial array  $\left\{ n \atop k\right\}_{F_n(x,s,q)} $ is not considered in  \cite{Cigler q-Lucas 2010}. Similarily:
\vspace{0.1cm}

\noindent \textbf{2.} the $q$-Lucas polynomial $L_n(x,s,q) = \sum_{c\in \Lambda_n}w(c) \equiv w(\Lambda)$ is the weight function of the set $\Lambda_n$ of  all coverings $c$ with arc monominos and dominos of the circle whose circumference has length $n$. Hence $L_n(x,s,q)$ is corresponding generation function for number of tilings of the circle whose circumference has length $n$. It may be then  combinatorially seen that $w(\Lambda_n)= w(\Phi_{n+1}) + s\cdot w(\Phi_{n-1}))$  hence $L_n(x,s,q) = F_{n+1}(x,s,q)+ s\cdot F_{n-1}(x,s,q)$. 

\noindent Polynomials $L_n(x,s,q)$ satify this paper recursion (3)  with $H_0(x)=2$ , $H_1(x)=x$ ; $s(x)=x$ and $t(x)=s$.

\noindent The $L_n(x,s,q)$-binomial array  $\left\{ n \atop k\right\}_{L_n(x,s,q)} $ is not considered in  \cite{Cigler q-Lucas 2010}.


\vspace{0.3cm}

\noindent Listing. \textbf{12.} In \cite[2010]{BSCS}   by  Bruce E. Sagan  and Carla  D.  Savage  the symbol $\left\{n\right\}\equiv U_n$    denotes the $n-th$ element of the  fundamental Lucas sequence $U$ satisfying this paper recurrence (2) with initial conditions $\left\{0\right\}=0$, $\left\{1\right\}=1$.  Naturally  $\left\{n\right\}$ is a polynomial in parameters  $s,t$. So is also the $U$-binomial coefficient $\left\{ n \atop k\right\}_U \equiv \left\{ n \atop k\right\}_{p,q}$. 

\vspace{0.1cm}

\noindent Similarly -  the symbol $\left\langle n \right\rangle\equiv V_n$    denotes the $n-th$ element of the  primordial Lucas sequence $V$ satisfying this paper recurrence (2) with initial conditions $\left\langle 0 \right\rangle = 2$ , $\left\langle 1 \right\rangle = s$ .  Naturally  $\left\langle n \right\rangle $  is a polynomial in parameters  $s,t$. So is also the $V$-binomial coefficient $\left\{ n \atop k\right\}_V \equiv \left\langle  n \atop k\right\rangle_{p,q}$.  $V$-binomials are not considered in \cite[2010]{BSCS}. Both fundamental and primordial sequences are interpreted via tilings similarly to the above in  \textbf{11.}  Johann Cigler attitude  rooted in already text-books tradition. 

\vspace{0.1cm}

\noindent An so:  $\left\{n\right\}$ is generation function for number of linear tilings  of  $(n-1)\times 1$  rectangle or equivalently of number of Morse code sequences of length $n-1$.

\vspace{0.1cm}

\noindent $\left\langle n \right\rangle $ is generation function for number of circular tilings of the circle whose circumference has length $n$. Using naturally proved (just seen) relations  Bruce E. Sagan  and Carla  D.  Savage  derive \textbf{two} combinatorial interpretations of the the same $\left\{ m+n \atop m,n\right\}_{p,q}$  via  Theorem 3.1. from which we infer the following.

\vspace{0.1cm}

\noindent \textbf{ 1.}   $\left\{ m+n \atop m,n\right\}_{p,q}$ is  the weight of all linear tilings of all integer partitions  $\lambda$ inside the $m \cdot n$ rectangle 

 hence $\left\{ m+n \atop m,n\right\}_{p,q}$ is the generating function  for numbers of such tilings of partitions.

\vspace{0.1cm}

\noindent  \textbf{2.}  $ 2^{m+n}\cdot \left\{ m+n \atop m,n\right\}_{p,q}$  is  the weight of all circular tilings of all integer partitions  $\lambda$ inside the $m \cdot n$ rectangle 

 hence $\left\{ m+n \atop m,n\right\}_{p,q}$ is the generating function  for numbers of such tilings of partitions.

\vspace{0.1cm}

\noindent \textbf{Explanation.} from \cite[2010]{BSCS}. \textit{A linear tiling of a partition}  $\lambda$   is a covering of its Ferrers diagram with disjoint dominos and monominos obtained by linearly tiling each  $\lambda_i$ part.  In circular tiling of a partition  $\lambda$ one performs circular tiling of each  $\lambda_i$ part

\vspace{0.3cm}

\noindent \textbf{The above list is open and far from complete.}


\vspace{0.5cm}

\noindent \textbf{4.2.} 

\noindent Nevertheless, to this end let us discern  in part- via indicative information - a part of Arthur T. Benjamin's recent contribution to the domain  . 
Namely; in  \cite[2003]{BQ2003}   by Arthur T. Benjamin and  Jennifer J. Quinn  track  the tilings' Combinatorial Theorem 5, p.36. There  for  $H_n = U_n$  the number   $s$ from the recurrence (2)  is interpreted as  equal to the number of colors of squares and  $t$ from this very recurrence (2)  equals to the number of colors of  dominos  while  $H_n = U_{n+1}$    counts colored tilings of length $n$  with squares and dominos.  
\noindent Similarly,  also in  \cite[2003]{BQ2003}  see   the  tilings'  Combinatorial Theorem 6 , p.36.  Here  for $H_n = V_n$  the number   $s$ from the recurrence (2)  should equal to the number of colors of a square and  $t$ from this very recurrence (2)  equals to the number of colors of a domino while  $H_n = V_n$  counts colored bracelets of length $n$  tiled with squares and dominos. 
\noindent  Bruce E. Sagan  and Carla  D.  Savage in  \cite[2010]{BSCS} refer to well known recurrences:  Identity 73 on p. 38 in \cite{BQ2003} - for (4) in \cite{BSCS} and Identity 94 p. 46 in \cite{BQ2003}  for (5) in \cite{BSCS}. Both (4) and (5) recurrences in \cite[2010]{BSCS} by Bruce E. Sagan  and Carla  D.  Savage  have been evoked  in the illustrative Example 3. Section 3. above.

\vspace{0.1cm}

\noindent Partially based on  \cite[2003]{BQ2003}  by Arthur T. Benjamin and Jennifer J. Quinn the paper \cite[2009]{BP}   by Arthur T. Benjamin and  Sean S. Plott should  be notified and as being nominated  by Arthur T. Benjamin and  Sean S. Plott in  errata \cite[2010]{BP} the present author feels entitled to remark on this errata.

\vspace{0.2cm}

\noindent \textbf{4.2.} According to errata \cite[2010]{BP} by  Arthur T. Benjamin and  Sean S. Plott [quote] " \textsl{The formula for}  ${n \choose k}_F$  \textit{should be multiplied by a factor of} $F_{n - x_k}$, \textsl{which accounts for the one remaining tiling that follows the $f_0$ tiling. Likewise, the formula for  ${n \choose k}_F$ should be multiplied by} $U_{ n - x_k}$."  Our remark is that this errata is unsuccessful. If we follow this errata  then ($x_{k-1} < x_k$) we would have: 

\begin{equation}
	\fnomialF{n}{k}{errata} =
	\sum_{1\leq x_1<x_2<\cdots<x_{k-1}\leq n-1}
	\prod_{i=1}^{k-1} 
	F_{k-i}^{x_i - x_{i-1}-1}
	F_{n-x_i - (k-i) + 1}
	F_{n - x_k},
\end{equation}
where   $F_0 = 0$ and $x_0 = 0$. But the formula (44) implies  for example

$$ 15 = {5 \choose 3}_F \neq \fnomialF{5}{3}{errata} = 11.$$

\vspace{0.1cm}

\noindent The task of finding the correct formula - due to the present author became a month ago an errand - exercise for Maciej Dziemia\'nczuk, a doctoral student from Gda\'nsk University in Poland.
The result - to be quoted below as MD formula (46)  - is his discovery, first announced in the form of a feedback private communication  to the present author: (M. Dziemia\'nczuk  on Mon, Oct 18, 2010 at 6:26 PM)  however still not announced in public.  

\vspace{0.1cm}

\noindent The source of an error in errata is that  ${n \choose k}_F$  should be multiplied \textbf{not by} the factor of $F_{n - x_k}$  \textit{but by} the factor  $F_{n - x_k +1}\equiv f_{n - x_k}$. Then we have

$$
	\fnomialF{n}{k}{now} =
	\sum_{1\leq x_1<x_2<\cdots<x_{k-1}\leq n-1}
	\prod_{i=1}^{k-1} 
	F_{k-i}^{x_i - x_{i-1}-1}
	F_{n-x_i - (k-i) + 1}
	F_{n - x_k +1},
$$
Due to $x_{k-1} < x_k$  the above  formula is equivalent to 

\begin{equation}
	\fnomialF{n}{k}{now} =
	\sum_{1\leq x_1<x_2<\cdots<x_{k-1<x_k}\leq n}
	\prod_{i=1}^{k-1} 
	F_{k-i}^{x_i - x_{i-1}-1}
	F_{n-x_i - (k-i) + 1}
	F_{n - x_k},
\end{equation}
\vspace{0.1cm}

\noindent and this in turn is evidently equivalent to the MD-formula (46) below i.e. (45) is equivalent to the corrected by Maciej Dziemia\'nczuk   Benjamin and Plott  formula from The Fibonacci Quarterly $46/47.1$ (2008/2009), 7-9. 

\vspace{0.2cm}

\noindent Finally here now  MD-formula follows:

\begin{equation}
	\fnomial{n}{k} =
	\sum_{1\leq x_1<x_2<\cdots<x_k\leq n}
	\prod_{i=1}^k 
	F_{k-i}^{x_i - x_{i-1}-1}
	F_{n-x_i - (k-i) + 1},
\end{equation}
where $F_0 = 0$ and $x_0 = 0$. 

\vspace{0.2cm}

\noindent Collaterally  Maciej Dziemia\'nczuk  supplies correspondingly correct formula for  Lucas $U$ - binomial  coefficients $\fnomialF{n}{k}{U}$ :

\begin{align}
	\fnomialF{n}{k}{U} &=
	\sum_{{1\leq x_1<x_2<\cdots<x_{k-1}\leq n} \atop {x_k = x_{k-1}+1}}
	\!\!\!\!\!\!\!\!\!\!\!\!
	s^{x_k - k}
	\left(\prod_{i=1}^{k-1} 
	U_{k-i}^{x_i - x_{i-1}-1}
	U_{n-x_i - (k-i) + 1}\right)
U_{n - x_{k} + 1}
	\\
	&= 
	\sum_{1\leq x_1<x_2<\cdots<x_k\leq n}
	\!\!\!\!\!\!\!\!
s^{x_k - k}
	\prod_{i=1}^k 
	U_{k-i}^{x_i - x_{i-1}-1}
	U_{n-x_i - (k-i) + 1},
\end{align}
where $U_0^t = 0^t = \delta_{t,0}$.


\vspace{0.6cm}

\noindent \textbf{4.3.} 
\noindent \textbf{$p,q$-binomials versus $q*$-binomials combinatorial interpretation,} where   $q* = \frac{p}{q}$ if  $q\neq 0$.

\vspace{0.2cm}

\noindent In the first instance let us  once for all switch off the uninspired  $p\cdot q = 0$ case. Then obligatorily either $q\neq 0$ or  $q\neq 0$. Let then  $q* = \frac{p}{q}$. In this nontrivial case

\begin{equation}
\fnomialF{n}{k}{p,q} =  q^{k(n-k)}\cdot\fnomialF{n}{k}{q*}.
\end{equation}
\vspace{0.2cm}

\noindent Referring to the factor $ q^{k(n-k)}$ as a kind of weight,  one may transfer combinatorial interpretation statements on $q*$ binomials $\fnomialF{n}{k}{q*}$ onto 
combinatorial interpretation statements on $p,q$ binomials $\fnomialF{n}{k}{p,q}$ through the agency of (49). Thence , apart from specific  combinatorial interpretations uncovered for the class or subclasses of $p,q$-binomials there might be admitted and respected the "$q*$-overall" combinatorial interpretations transfered from $1,q*$-binomials i.e. from $q*$-binomials onto $p,q$-binomials.

\vspace{0.2cm}

\noindent By no means pretending to be the complete list here comes the skeletonized list of \textbf{[Ex. q* ; k]} examples,  $k\geq 1$.

\vspace{0.4cm}

\noindent \textbf{[Ex. q* ; 1]}  
\noindent The $q*$-binomial coefficient ${m+n \choose m,n}_{q*}$  may be interpreted as a polynomial in $q*$ whose $q*^k$-th coefficient   counts the number of distinct partitions of $k$ elements which fit inside an $m \times n$ rectangle - see  \cite[1976]{Andrews}  by   George Eyre Andrews.

\vspace{0.3cm}

\noindent \textbf{On lattice path techniques  - Historical Remark.} It seams to  be desirable  now to  quote here information from \cite[2010]{Katherine 2010}  by Katherine Humphreys based on
\cite[1878]{Whitworth 1878}  by William Allen Whitworth:

\begin{quot}
We find lattice path techniques as early as 1878 in Whitworth to help picture a combinatorial problem, but it is not until the early 1960's that we find lattice path enumeration presented as a mathematical topic on its own. The number of papers pertaining to lattice path enumeration has more than doubled each decade since 1960. 
\end{quot}

\vspace{0.3cm}

\noindent \textbf{[Ex. q* ; 2]} 
\noindent The [Ex. q* ; 2] may be now compiled  with [Ex. q* ; 1] above. For that to do recall that zigzag path is the shortest path  that  starts at
$A = (0,0)$  and ends in $B =  (n,n-k)$ of the $n \times k$ rectangle; see:  \cite[1962]{P62}  by  Gy{\"{o}}rgy P\'olya   [pp.  68-75],   \cite[1969]{P62}  by  Gy{\"{o}}rgy P\'olya  and
\cite{PA 1971} by  Gy{\"{o}}rgy P\'olya    and  G. L.  Alexanderson.

\vspace{0.1cm}

\noindent Let then     $A_{n,k,\alpha}= $ the number of those $(0,0) \longrightarrow  (k,n-k)$ zigzag paths the area under which  is $\alpha$.

\vspace{0.1cm}

\noindent In  \cite[1969]{P69}  Gy{\"{o}}rgy  P\'olya  using recursion for  $q*$-binomial coefficients proved  that  

$$  \fnomialF{n}{k}{q*} = \sum_{\alpha = 0}^{k(n-k)} A_{n,k,\alpha}\cdot q*^{\alpha}.$$    
from where Gy{\"{o}}rgy  P\'olya infers the following Lemma (\cite[1969]{P69}, p.105) which is named Theorem (p. 104) in more detailed paper \cite[1971]{PA 1971} by  Gy{\"{o}}rgy  P\'olya  and  G. L.  Alexanderson.

\begin{quot}
The number of those zigzag paths the area under which  is $\alpha$  equals  $A_{n,k,\alpha}$.
\end{quot}

\vspace{0.3cm}

\noindent \textbf{[Ex. q* ; 3]} 
\noindent The [Ex. q* ; 3] may be now compared  with [Ex. q* ; 1]. The combinatorial interpretation of ${r+s \choose r,s}_{q*}$ from  [Ex. q* ; 1] had been derived (pp. 106-107) in \cite[1971]{PA 1971} by  Gy{\"{o}}rgy  P\'olya  and  G. L.  Alexanderson, from where - with advocacy from \cite[1971]{Knuth1971}   by  Donald  Ervin  Knuth - we quote the result.

\vspace{0.4cm}

\noindent \textbf{(1971)}: \; ${r+s\choose r,s}_{q*}=$ \textit{ordinary generating function in $\alpha$  powers  of  $q*$  for partitions of $\alpha$ into exactly $r$  non-negative integers none of which exceeds $s$} ,  
\vspace{0.2cm}

\noindent as derived  in \cite[1971]{PA 1971} by  Gy{\"{o}}rgy  P\'olya  and  G. L.  Alexanderson  - see formula (6.9) in \cite{PA 1971}.

\vspace{0.4cm}

\noindent \textbf{(1882)}:\; $\fnomialF{n}{k}{q*} = $ \textit{ordinary generating function in $\alpha$  powers  of  $q$  for partitions of $\alpha$ into at most $k$ parts not exceeding $(n-k)$} ,  

\vspace{0.2cm}

\noindent as recalled in \cite[1971]{Knuth1971} by  Donald  Ervin  Knuth and proved combinatorially in \cite[1882]{SJJ 1882} by James Joseph Sylvester.

\vspace{0.2cm}

\noindent  Let nonce : \: $r+s=n$, $r=k$ then $\textbf{(1971)}\equiv \textbf{(1882)}$ are equal due to    

\begin{equation}
\fnomialF{n}{k}{q*} = \sum_{\alpha = 0}^{k(n-k)} A_{n,k,\alpha}\cdot q*^{\alpha} = \sum_{\alpha = 0}^{r\cdot s} A_{r+s,r,\alpha}\cdot q*^{\alpha} = {r+s \choose r,s}_{q*}.
\end{equation}
where for commodity of comparison formulas in two notations from two papers -  we have been using contractually for a while: $r+s=n$, $r=k$ identifications.

\vspace{0.3cm}

\noindent \textbf{[Ex. q* ; 4]}  
\noindent The following was proved in \cite[1961]{Kendall 1961} by Maurice  George Kendall and  Alan Stuart (see p.479 and p.964) and n \cite[1971]{PA 1971} by  Gy{\"{o}}rgy  P\'olya  and  G. L.  Alexanderson  (p.106).

\vspace{0.3cm}

\noindent \textit{The area under the zigzag path $=$  The number of inversions  in  the very zigzag path coding sequence.}

\vspace{0.2cm}

\noindent The possible extension of the above combinatorial interpretation onto three dimensional zigzag paths via  "\textit{three-nomials}" was briefly mentioned in \cite{PA 1971} - see p.108.

\vspace{0.3cm}

\noindent \textbf{[Ex. q* ; 5]}  
\noindent The well known (in consequence - \textit{finite geometries'}) interpretation of $\fnomialF{n}{k}{q*}$ coefficient due to Jay Goldman and  Gian-Carlo Rota from \cite[1970]{Rota Goldman 1970} is now  worthy of being recalled;  see also  \cite[1971]{Knuth1971} by  Donald  Ervin  Knuth.

\vspace{0.1cm}

\noindent Let   $V_n$ be an $n$-dimensional vector space over a finite field of $q*$ elements. Then 

\vspace{0.2cm}

\noindent $\fnomialF{n}{k}{q*} $ = \textit{the number of $k$-dimensional subspaces of  $V_n$ }.

\vspace{0.3cm}

\noindent \textbf{[Ex. q* ; 6]}  

\noindent This example $=$ the short substantial note  \cite[1971]{Knuth1971} by  Donald  Ervin  Knuth. 

\noindent Compile this example with the example [Ex. q* ; 5] above.

\noindent The essence of a coding  of combinatorial interpretations via bijection between lattices  is the construction of this \textit{coding bijection}  in \cite{Knuth1971}. Namely, let $GF(q*)$ be the Galois field of order $q*$ and let $ V_n \equiv V = GF(q*)^n$ be the $n$-dimensional vector space over $GF(q*)$. Let $[n]= \left\{1,2,...,n\right\}$. Let $\ell(V)$ be the lattice of all subspaces of $ V = GF(q*)^n$  while  $\ell([n])\equiv 2^{[n]}$ denotes the lattice of all subsets of $ [n]$.

\noindent In \cite{Knuth1971} Donald  Ervin  Knuth constructs this  \textit{natural order and rank preserving} map $\Phi$  from the lattice $\ell (V)$ of subspaces  onto the lattice $\ell ([n])\equiv 2^{[n]}$  of subsets of $[n]$. 
  
$$  \ell(V)\;   \stackrel{\Phi}{\rightarrow}\; \ell([n]).$$
\noindent We bethink with some reason whether  this $\Phi$ bijection coding  might be an answer to the subset-subspace problem from subset-subspace problem  from \cite[1998]{Konva 1998} by John Konvalina \textbf{?}

\begin{quot}
 ...the subset-subspace problem (see 6 , 9 , and 3) . The traditional approach to the subset-subspace problem has been to draw the following analogy:   the binomial $\fnomial{n}{k}$  coefficient  counts k-subsets of an n-set, while the analogous Gaussian $\fnomialF{n}{k}{q}$  coefficient counts the number of $k$-dimensional subspaces of an $n$-dimensional finite vector space over the field of 
$q$ elements. 
\noindent The implication from this analogy is that the Gaussian coefficients and related identities tend to the analogous identities for the ordinary binomial coefficients as $q$ approaches 1.
The proofs are often algebraic or mimic subset proofs. But what is the combinatorial reason for the striking parallels between the Gaussian coefficients and the binomial coefficients?
\end{quot}

\noindent According to  Joshef  P. S. Kung  \cite[1995]{Kung 1995} the Knuth's note is not the explanation: 

\begin{quot}
... observation of Knuth yields an order preserving map from $L(V_n(q)$ to Boolean algebra of subsets, but it does not yield a solution to the still unresolved problem of finding a combinatorial interpretation of taking the limit $q \longrightarrow 1$.
\end{quot}Well, perhaps  this limit being performed by $q$-deformed Quantum Mechanics  physicists might be of some help?  There the so called $q$-quantum plain of $q$-commuting variables $x\cdot y- q\cdot y\cdot x =0$   becomes a plane $\mathbb{F}\times \mathbb{F}$  ($\mathbb{F}\  =  \mathbb{R}, \mathbb{C} $,...  $p$-adic fields included) of two commuting variables in the limit $q \longrightarrow 1$. For   see \cite[1953]{ Schutzenberger 1953}  by Marcel-Paul   Sch{\"{u}}tzenberger. For quantum plains - see also  \cite[1995]{Kassel 1995}  by  Christian Kassel. It may deserve notifying that $q$ - extension of of the "classical plane" of commuting variables ($q=1$) seems in a sense ultimate as discussed in \cite[2001]{akk 2001} by  A.K. Kwa\'sniewski

\vspace{0.3cm}

\noindent \textbf{[Ex. q* ; 7]}  
\noindent Let us continue the above by further quotation from \cite[1998]{Konva 1998} on generalized binomial coefficients and the subset-subspace problem.

\begin{quot}  
We will show that interpreting the Gaussian coefficients as generalized binomial coefficients of the second kind combinations with repetition reveals the combinatorial connections between not only the binomial coefficients and the Gaussian coefficients, but the Stirling numbers as well. Thus, the ordinary Gaussian coefficient  tends to be an algebraic generalization of the binomial coefficient of the first kind, and a combinatorial generalization of the binomial coefficient of the
second kind.
\end{quot}
Now in order to get more oriented go back to the begining  of  subsection \textbf{ 4.1.} and consult  : Listing. \textbf{1.}, Listing. \textbf{2.}, Listing. \textbf{3.}  which are earlier works and
end up with  \cite[2000]{Konva 2000}  by   John Konvalina on an unified  simultaneous  interpretation of binomial coefficients of both kinds, Stirling numbers of both kinds and Gaussian binomial coefficients  of both kinds.  Compare it then afterwards  with  Listing. \textbf{8.}.

\vspace{0.3cm}




\end{document}